\documentclass[12pt,english,refpage,intoc,bibliography=totoc,index=totoc,BCOR7.5mm,captions=tableheading]{paper}
\usepackage[T1]{fontenc}
\usepackage[latin9]{inputenc}
\usepackage{color}
\usepackage{babel}
\usepackage[unicode=true,
 bookmarks=true,bookmarksnumbered=true,bookmarksopen=false,
 breaklinks=false,pdfborder={0 0 1},backref=false,colorlinks=true]
 {hyperref}
\hypersetup{pdftitle={The LyX User's Guide},
 pdfauthor={LyX Team},
 pdfsubject={LyX},
 pdfkeywords={LyX},
 linkcolor=black, citecolor=black, urlcolor=blue, filecolor=blue, pdfpagelayout=OneColumn, pdfnewwindow=true, pdfstartview=XYZ, plainpages=false}
\usepackage{breakurl}

\makeatletter
\usepackage{enumitem}		

\usepackage{ifpdf} 
\ifpdf 

 \IfFileExists{lmodern.sty}{\usepackage{lmodern}}{}

\fi 

\usepackage[figure]{hypcap}

\let\myTOC\tableofcontents
\renewcommand\tableofcontents{%
  \frontmatter
  \pdfbookmark[1]{\contentsname}{}
  \myTOC
  \mainmatter }

\usepackage{fancyhdr}

\makeatother

\begin{document}

\title{A Proof of the Isometric Embedding Theorem in Euclidean Space of
Dimension Three.}

\author{Edgar Kann}

\date{November 21, 2017}
\maketitle
\begin{abstract}
A proof of the isometric embedding of a given two-metric in $E^{3}$
of low differentiability class. The method uses the theory of first
order partial differential equations. The curvature of the metric
plays no role in the proof. 
\end{abstract}

\section*{Introduction.}

Non-linear first order systems of partial differential equations giving
the transformation of the components of a metric expressed in one
parameter system to its expression in components in another parameter
system play a vital role in the embedding problem. Although the systems
are non-linear in the partial derivatives they are linear in the algebraic
sense in the components of the metric. Thus we can make use of the
consistency theorem of linear algebra, namely: a linear algebraic
system is consistent (i.e. has a solution) if and only if the rank
of the augmented matrix is equal to the rank of the coefficient matrix.
In our case the elements of the coefficient matrix are the squares
of the partial derivatives of the parameters of one parameter system
with respect to parameters in a second parameter system. The consistency
theorem thus gives rise to two linear partial differential equations
of first order the existence of whose solutions is sufficient for
the existence of solutions of the original non-linear PDE systems.
The applications of this idea are given in two parts in this paper.
In both parts the metric discussed is presented in orthogonal parameters.

In Part Two we give sufficient conditions on the components of any
presented metric for it to be isometrically embedded as a surface
in $E^{3}.$ The sufficient conditions are in the form of two equations
involving the solutions of two linear PDEs of first order whose solutions
are determined by respective initial conditions. The initial conditions
may be any $C^{1}$ functions and are treated as unknown functions
to be determined using Parts Two and Three.

In Part Three we present sufficient conditions on functions that they
be components of a given metric presented $a'priori$ in geodesic
parameters. This part is purely intrinsic.

In Part Four we put the two parts together in two differential equations
in which the initial conditions of Part Two play the role of unknown
functions. The two differential equations in Part Four form a system
of ordinary differential equations and so can be solved by appropriate
initial conditions. Because the sufficiency condition of Part Two
is satisfied by the system some isometric embedding exists. Because
the sufficiency condition of Part Three is satisfied as well this
embedding is an isometric embedding of the $\mathit{a'priori}$ given
metric. 

Part Five presents a simple example illustrating the method.

\section{Main theorem.}

Any 2-metric can be locally isometrically embedded in $E^{3}$ in
the form 
\[
X(u,v)=x(u,v)\mathbf{i}+y(u,v)\mathbf{j}+v\mathbf{k}
\]
 where $u$ and $v$ are orthogonal parameters on the embedded surface
and $\mathbf{i},\mathbf{j}$,$\mathbf{k}$ is an orthonormal basis
of $E^{3}.$

We assume that the components of the metric are in differentiability
class $C^{1}.$ We refer to $u$ and $v$ as $\mathit{level}$ $\mathit{parameters.}$

\section{Preliminary. Isometric Embedding of a special metric. }

Consider the system $S_{0}$ defined as 
\[
E=x_{u}^{2}+y_{u}^{2}
\]
\[
0=x_{u}x_{v}+y_{u}y_{v}
\]
\[
G-1=x_{v}^{2}+y_{v}^{2}
\]
 with $G>1$, $E$ and $G$ functions of $u,v.$ 

Theorem 1. If $S_{0}$ is satisfied by $C^{1}$ functions $x(u,v),y(u,v),E(u,v),G(u,v),$
then the surface defined above by $X(u.v)$ is an isometric embedding
of the metric $\omega=E(u,v)du^{2}+G(u,v)dv^{2}.$ 

Proof. The metric components induced by the ambient space $E^{3}$
on $X$ are:

\[
X_{u}^{2}=x_{u}^{2}+y_{u}^{2}
\]
\[
X_{u}\cdot X_{v}=x_{u}x_{v}+y_{u}y_{v}
\]
\[
X_{v}^{2}=x_{v}^{2}+y_{v}^{2}+1.
\]
 Therefore $E=X_{u}^{2},\,$$F=X_{u}\cdot X_{v}=0,\, G=X_{v}^{2}.$ 

Thus $X$ is, in fact, a regular surface since $EG-F^{2}>0.$ Hence,
$X$ is an isometric embedding of the metric $\omega$. QED theorem
1. 

The system $S_{0}$ says that the Euclidean plane metric $\omega_{0}=dx^{2}+dy^{2}$
has components $E(u,v)$, 0, $G(u,v)-1$ in parameter system $u,v$
and the functions $x(u,v),\, y(u,v)$ give the transformation of parameters
from the u,v parameters to the canonical $x,y$ parameters of the
plane. 

Consider the inverse system $S_{0}^{-1}$ which says that the plane
metric $\omega_{0}$ has components $1,0,1$ in canonical parameters
$x,y$ and the functions $u(x,y),\, v(x,y)$ give the transformation
of parameters from the $x,y$ parameters to the $u,v$ parameters.
$S_{0}^{-1}$ is the system 
\begin{equation}
1=E(u,v)u_{x}^{2}+(G(u,v)-1)v_{x}^{2}.\label{S_01}
\end{equation}
\begin{equation}
0=E(u,v)u_{x}u_{y}+(G(u,v)-1)v_{x}v_{y}.\label{S_02}
\end{equation}
\begin{equation}
1=E(u,v)u_{y}^{2}+(G(u,v)-1)v_{y}^{2}.\label{S_03}
\end{equation}
 $S_{0}$ has a solution $x(u,v),y(u,v),E(u,v),G(u,v),$ if and only
if $S_{0}^{-1}$ has a solution $E(u,v),\, G(u.v),\, u(x,y),\, v(x,y).$
The functions $x(u,v),\, y(u,v)$ give the transformation of parameters
from the u,v parameters to the canonical $x,y$ parameters of the
plane. The functions $u(x,y),\, v(x,y)$ give the inverse parameter
transformation. We will work with $S_{0}^{-1}$ to show that it has,
in fact, such solutions. We make use of the theory of first order
partial differential equations. \cite{John}

\label{Lemma-1.-The}Lemma 1. The two linear first order PDEs
\begin{equation}
u_{x}-u_{y}=0,\label{eq:first order PD for system S0 inverse-3}
\end{equation}
\begin{equation}
v_{x}+v_{y}=0.\label{eq:first order PDE for system S0 inverse}
\end{equation}
 have Cauchy problem solutions 
\begin{equation}
u=a(x+y),\label{eq:sol 1st Cauchy prob}
\end{equation}
\begin{equation}
v=b(x-y)\label{Sol 2nd Cauchy prob}
\end{equation}

with initial conditions respectively on $y=0:$ $u=a(x)$ and $v=b(x).$
$a$ and $b$ are arbitrary $C^{1}$ functions chosen to have non-zero
first derivatives. The proof is by inspection. QED Lemma 1.

The two solutions provide a one-to-one parameter transformation from
the $xy$ plane to the $uv$ plane whose Jacobian evaluated on $y=0$
has the value
\[
J(x,0)=-2a'(x)b'(x),
\]
 where the prime indicates the first derivative. The parameter transformation
is thus locally invertible in a neighborhood of the initial curve.

Noting that $x$ and $y$ are functions of $u$ and $v$, we define

\begin{equation}
E(u,v)=\frac{1}{2}(a'(x+y))^{-2}.\label{eq:Def E}
\end{equation}
\begin{equation}
G(u,v)-1=\frac{1}{2}(b'(x-y))^{-2}.\label{Def G}
\end{equation}
 The data space of functions $a(x),b(x)$ is thus mapped to the space
of functions $E(u,v),G(u,v).$

Claim. The system $S_{0}^{-1}$ is satisfied by the functions $E(u,v),\, G(u.v),\, u(x,y),\, v(x,y).$

Proof of Claim. By substitution. See subsection 2.1 below on calculation.
QED Claim.

Therefore $S_{0}$ has a solution $x(u,v),y(u,v),E(u,v),G(u,v).$ 

By Theorem 1

\[
X(u,v)=x(u,v)\mathbf{i}+y(u,v)\mathbf{j}+v\mathbf{k}
\]

is an isometric embedding in $E^{3}$ of the metric defined by
\[
\omega=Edu^{2}+Gdv^{2}.
\]

\subsection{Calculation.}

Substitute from the definitions of $E,$ $G$ into the right sides
of the equations of systems $S_{0}^{-1}$ using the solutions ( \ref{eq:sol 1st Cauchy prob},
\ref{Sol 2nd Cauchy prob} ) of the PDEs to show that the equations
are satisfied.
\[
\frac{1}{2}(a'(x+y))^{-2}u_{x}^{2}+(G(u,v)-1)v_{x}^{2}
\]
\[
=\frac{1}{2}(a'(x+y))^{-2}(a'(x+y))^{2}+\frac{1}{2}(b'(x-y)^{-2}(b'(x-y)^{2}=1.
\]
 Therefore the first equation of $S_{0}^{-1}$ is satisfied.

\[
\frac{1}{2}(a'(x+y))^{-2}a'(x+y)a'(x+y)+\frac{1}{2}(b(x-y)^{-2}b'(x-y)b'(x-y)(-1)=0.
\]
 Therefore the second equation of $S_{0}^{-1}$ is satisfied. 
\[
\frac{1}{2}(a'(x+y)^{-2}(a'(x+y)^{2}+\frac{1}{2}(b'(x-y)^{-2}(b'(x-y)(-1))^{2}=1.
\]
 Therefore the third equation of $S_{0}^{-1}$ is satisfied. QED Claim.

The parameter transformation (6),$\,$(7) has non-zero Jacobian in
a neighborhood of the initial curve and thus can be inverted there.
Denote the inverse parameter transformation by
\[
x=H(u,v),\, y=H^{*}(u,v).
\]
 Thus by Theorem 1, we have proved 

Theorem 2. \label{Theorem-2.-(}

\[
X(u,v)=H(u,v)\mathbf{i}+H^{*}(u,v)\mathbf{j}+v\mathbf{k}
\]
 is an isometric embedding in $E^{3}$ of the metric 
\[
\omega=Edu^{2}+Gdv^{2}
\]
 where $E$ and $G$ are defined above depending on the choice of
initial conditions $a$ and $b.$ (That is, $E$ and G are not given
functions a' priori.)

\section{Transformation of components of a given metric expressed in geodesic
parameters.\label{sec:Transformation-of-components}}

We assume that the components of the given metric are in class $C^{1}.$
We define the given metric by
\[
\omega=d\hat{u}^{2}+\hat{G}(\hat{u},\hat{v})d\hat{v}^{2}.
\]
Transform the metric to a system of orthogonal parameters $u,v$ in
which the components of $\omega$ are $R(u,v),\,0,\, S(u,v):$ 
\[
\omega=Rdu^{2}+Sdv^{2}.
\]
 Then the components in the two systems are related by a system $\hat{S}$
defined as: 
\[
1=Ru_{\hat{u}}^{2}+Sv_{\hat{u}}^{2}.
\]
\[
0=Ru_{\hat{u}}u_{\hat{v}}+Sv_{\hat{u}}v_{\hat{v}}.
\]
\[
\hat{G}(\hat{u},\hat{v})=Ru_{\hat{v}}^{2}+Sv_{\hat{v}}^{2}.
\]
 The system is algebraically linear in the unknowns $R,S$ with the
squared derivatives considerd as coefficients to be determined independently
below (\ref{eq:linear PDE for u.}, \ref{linear PDE for v.}). Therefore
a necessary and sufficient condition for there to be an algebraic
solution for $R$,$\mathrm{S}$ is that the rank of the augmented
matrix equals the rank of the coefficient matrix. The augmented matrix
is
\[
\left[\begin{array}{ccc}
u_{\hat{u}}^{2} & v_{\hat{u}}^{2} & 1\\
u_{\hat{u}}u_{\hat{v}} & v_{\hat{u}}v_{\hat{v}} & 0\\
u_{\hat{v}}^{2} & v_{\hat{v}}^{2} & \hat{G}(\hat{u},\hat{v})
\end{array}\right].
\]

The determinant of the augmented matrix is
\[
\left|\begin{array}{ccc}
u_{\hat{u}}^{2} & v_{\hat{u}}^{2} & 1\\
u_{\hat{u}}u_{\hat{v}} & v_{\hat{u}}v_{\hat{v}} & 0\\
u_{\hat{v}}^{2} & v_{\hat{v}}^{2} & \hat{G}(\hat{u},\hat{v})
\end{array}\right|=\left|\begin{array}{cc}
u_{\hat{u}}u_{\hat{v}} & v_{\hat{u}}v_{\hat{v}}\\
u_{\hat{v}}^{2} & v_{\hat{v}}^{2}
\end{array}\right|+\begin{array}{c}
\hat{G}(\hat{u},\hat{v})\end{array}\left|\begin{array}{cc}
u_{\hat{u}}^{2} & v_{\hat{u}}^{2}\\
u_{\hat{u}}u_{\hat{v}} & v_{\hat{u}}v_{\hat{v}}
\end{array}\right|
\]
\[
=u_{\hat{v}}v_{\hat{v}}\left|\begin{array}{cc}
u_{\hat{u}} & v_{\hat{u}}\\
u_{\hat{v}} & v_{\hat{v}}
\end{array}\right|+\hat{G}(\hat{u},\hat{v})u_{\hat{u}}v_{\hat{u}}\left|\begin{array}{cc}
u_{\hat{u}} & v_{\hat{u}}\\
u_{\hat{v}} & v_{\hat{v}}
\end{array}\right|
\]
\[
=\left(u_{\hat{v}}v_{\hat{v}}+\hat{G}(\hat{u},\hat{v})u_{\hat{u}}v_{\hat{u}}\right)\left|\begin{array}{cc}
u_{\hat{u}} & v_{\hat{u}}\\
u_{\hat{v}} & v_{\hat{v}}
\end{array}\right|=0.
\]
We will insure (see Claim below) that the determinant is not zero. 

Thus the consistency condition is
\[
\left(u_{\hat{v}}v_{\hat{v}}+\hat{G}(\hat{u},\hat{v})u_{\hat{u}}v_{\hat{u}}\right)=0.
\]
 The partial derivatives are to be determined as solutions of Cauchy
problems for the first order linear PDEs
\begin{equation}
u_{\hat{u}}-u_{\hat{v}}=0,\label{eq:linear PDE for u.}
\end{equation}
\begin{equation}
\hat{G}(\hat{u},\hat{v})v_{\hat{u}}+v_{\hat{v}}=0.\label{linear PDE for v.}
\end{equation}
 Solutions of these two PDEs satisfy the consistency condition as
can be seen by substitution. Note, however, that they they are not
necessary for consistency: the identity transformation satisfies the
system $\hat{S}$ with $R=S=1$ but does not satisfy the PDEs. 

Solutions of Cauchy problems for these PDEs provide a solution of
the system $\hat{S}$ for $R$ and $S$ by algebraic methods: Substitute
the solutions into $\hat{S}$ and solve by Cramer's Rule for $R$
and $S.$ Since the system is consistent we need only two of the equations.
We choose the first two. 
\[
1=Ru_{\hat{u}}^{2}+Sv_{\hat{u}}^{2}.
\]
\[
0=Ru_{\hat{u}}u_{\hat{v}}+Sv_{\hat{u}}v_{\hat{v}}.
\]
 Substitute from the first order PDEs:
\[
1=Ru_{\hat{v}}^{2}+S\left(\frac{v_{\hat{v}}^{2}}{\hat{G}^{2}}\right).
\]
\[
0=Ru_{\hat{v}}^{2}-S\frac{v_{\hat{v}}^{2}}{\hat{G}}.
\]
Solve by Cramer's rule: 
\[
R(u,v)=u_{\hat{v}}^{-2}\frac{\hat{G}(\hat{u},\hat{v})}{\hat{G}(\hat{u},\hat{v})+1}.
\]
 
\[
S(u,v)=v_{\hat{v}}^{-2}\frac{\hat{G}^{2}(\hat{u},\hat{v})}{\hat{G}(\hat{u},\hat{v})+1}.
\]
 The theory of first order PDEs is explained clearly in \cite{John}.
The theory is particularly simple for the linear PDEs above. We solve
initial value problems for these. A general solution of the initial
value problem for the first PDE (\ref{eq:linear PDE for u.}) above
is worked out in the first problem in John's book, p. 15. With appropriate
change of notation the solution is 
\begin{equation}
u=h(\hat{u}+\hat{v})\label{eq:sol first PDE for u}
\end{equation}
with initial conditions on $\hat{v}=0:$ $u=h(\hat{u}).$ We stipulate
that the initial value satisfies $h_{\hat{u}}(\hat{u})\equiv h'(\hat{u})\neq0.$

A general solution of the second (\ref{linear PDE for v.}) is:
\begin{equation}
v=\hat{h}(\hat{u},\hat{v})\label{sol second PDE for v}
\end{equation}
 with chosen initial value $\hat{h}(\hat{u})=\hat{h}(\hat{u},0)$.
We also choose the initial value to satisfy that $v_{\hat{u}}(\hat{u},0)$
is not zero. 

Thus, for all $\hat{u},$ $\hat{v}$
\[
u_{\hat{v}}=h'(\hat{u}+\hat{v}),\, v_{\hat{v}}=\hat{h}_{\hat{v}}(\hat{u},\hat{v}).
\]
 Claim. The Jacobian of the transformation of parameters given by
the solutions of the two PDES is not zero at the initial curve (and
hence in a neighborhood of the initial curve). 

Proof of Claim. : The IC implies $\hat{h}_{\hat{v}}$ is not zero
(see \ref{linear PDE for v.}).

\[
J=\left|\begin{array}{cc}
u_{\hat{u}} & u_{\hat{v}}\\
v_{\hat{u}} & v_{\hat{v}}
\end{array}\right|=\left|\begin{array}{cc}
h'(\hat{u}) & h'(\hat{u})\\
\hat{h}_{\hat{u}} & \hat{h}_{\hat{v}}
\end{array}\right|=h'(\hat{u})(\hat{h}_{\hat{v}}-\hat{h}_{\hat{u}})=h'(\hat{u})(v_{\hat{v}}-v_{\hat{u}})=h'(\hat{u})[-\hat{G}(\hat{u},\hat{v})v_{\hat{u}}-v_{\hat{u}}].
\]
 Thus $J=h'(\hat{u})[-v_{\hat{u}}(\hat{G}+1)]\neq0.$ 

QED Claim.

Hence, using the solutions, 
\[
R(u,v)=[u_{\hat{v}}]^{-2}\frac{\hat{G}}{\hat{G}+1}.
\]
\[
S(u,v)=[\hat{h}_{\hat{v}}(\hat{u},\hat{v})]^{-2}\frac{\hat{G}^{2}}{\hat{G}+1}.
\]

Thus $R(u,v)$ is known because it is determined by the given metric
and the initial conditions. The value of $S(u,v)$ depends on the
choice of the IC on $\hat{v}$=0, that is, on the choice of $\hat{h}(\hat{u},0)$
such that $\hat{h}_{\hat{u}}(\hat{u},0)$ is not zero. 

Note that $R$ and $S$ are given here as functions of $\hat{u}$
and $\hat{v}$. We want them as functions of $u$ and $v.$ By what
we have just shown, the transformation of parameters is locally one-to-one
and has a $C^{1}$ inverse transformation giving $\hat{u},$ $\hat{v}$
in terms of $u,v$. 

Represent the inverse by:
\begin{equation}
\hat{u}=f(u,v),\label{eq:Def of f}
\end{equation}
\begin{equation}
\hat{v}=g(u,v)\label{eq:Def of g}
\end{equation}
 and substitute into the expressions for $R$ and $S:$ 
\[
R(u,v)=[u_{\hat{v}}]^{-2}\frac{\hat{G(}f(u,v),g(u,v))}{\hat{G(}f(u,v),g(u,v))+1}.
\]
\[
S(u,v)=[\hat{h}_{\hat{v}}(f(u,v),g(u,v)]^{-2}\frac{[\hat{G(}f(u,v),g(u,v))]^{2}}{\hat{G(}f(u,v),g(u,v))+1}.
\]
 Thus, $R(u,v)$ and $S(u,v)$ are determined by the choice of the
initial conditions $\hat{h}(\hat{u})=\hat{h}(\hat{u},0)$ and $h(\hat{u})=h(\hat{u}+0).$ 

Thus we have proved 

Theorem 3. $R(u,v)$ and $S(u,v)$ are components of $\omega.$ 

Remark. We could have simplified the argument a bit by choosing the
initial condition for (\ref{eq:linear PDE for u.}) such that the
solution is $u=\hat{u}+\hat{v}$ but that would have precluded the
change of metric components to other geodesic parameters (e.g., it
would produce a contradiction when $R=1).$ Since this has no consequences
for our embedding proof but simplifies the notation, in what follows
we will assume that the solution is $u=\hat{u}+\hat{v}.$

\section{Proof of the embedding theorem.}

We return to Lemma 1, section 2.

The parametric form of the solutions of the PDEs in Lemma 1 is as
follows.

For the first PDE (\ref{eq:first order PD for system S0 inverse-3})
we choose parameters $s=x+y,\: t=y.$

A parametric solution of the first PDE is (with John's notation $z$
replaced by $u$ in p. 15, example 1, \cite{John}):
\[
x=s-t,
\]
\[
y=t,
\]
\[
u=a(s).
\]
 The initial curve, on which $t=0,$ is
\[
x=s,
\]
\[
y=0,
\]
\[
u=a(s).
\]

For the second PDE (\ref{eq:first order PDE for system S0 inverse})
we choose parameters $\sigma=x-y,\:\tau=y.$ This is possible because
the two PDEs are independent. 

\[
x=\sigma+\tau,
\]

\[
y=\tau,
\]
\[
v=b(\sigma).
\]
 The initial curve for the second PDE, on which $\tau=0,$ is 
\[
x=\sigma,
\]
\[
y=0,
\]
\[
v=b(\sigma).
\]

(With $z$ replaced by $v$ in \cite{John}.)

Note that the characteristic curves of the PDEs satisfy respectively
$\frac{du}{dt}=0,\,\frac{dv}{d\tau}=0.$ Compare \cite{John}, p.14.
Also, note that the solution functions have the same form as the initial
conditions. Therefore we have 

Lemma 2. In results making use of solutions of the PDEs we may assume,
without loss of generality, that $y=0,\, x=s=\sigma,$ and $t=\tau=0.$

Consider the ODE system 
\[
\frac{1}{2}(a'(s))^{-2}=\frac{\hat{G(}f(a(s),b(s)),g(a(s),b(s)))}{\hat{G(}f(a(s),b(s)),g(a(s),b(s)))+1}.
\]

\[
\frac{1}{2}(b'(s))^{-2}+1=[\hat{h}_{\hat{v}}(f(a(s),b(s)),g(a(s),b(s)))]^{-2}\frac{[\hat{G(}f(a(s),b(s)),g(a(s),b(s)))]^{2}}{\hat{G(}f(a(s),b(s)),g(a(s),b(s)))+1}.
\]
 Put in standard form $a'=....,\, b'=,...$ choosing the plus signs
for the square roots, e.g. For given initial values of $a$ and $b$
at a point there exists a solution $a(s),\, b(s)$ of the system in
a neighborhood of the initial point satisfying the initial conditions.
\cite{Cod and Lev}

We choose the initial point to be $s=0$ and the initial values of
the unknown functions to be $a(0)=0,\, b(0)=0.$ Thus a solution of
exists satisfying these initial conditions.

Claim. The standard forms are real. $a'$ is real by inspection. For
$b'$ we must show that the right side of the second equation is greater
than one at the initial point. This is equivalent to showing that 

\[
[\hat{h}_{\hat{v}}(f(0,0),g(0,0))]^{2}<\frac{\hat{G}^{2}(f(0,0),g(0,0))}{\hat{G}^{2}(f(0,0),g(0,0))+1}.
\]
 By using the partial differential equations and their solutions from
the preceding section (\ref{eq:linear PDE for u.}, \ref{linear PDE for v.},
\ref{eq:sol first PDE for u}, \ref{sol second PDE for v}) we find
that

$\hat{h}_{\hat{v}}(\hat{u},0)=v_{\hat{v}}(\hat{u},0)=-\hat{G}(\hat{u},0)v_{\hat{u}}(\hat{u},0).$ 

By choosing the initial value $\hat{h}(\hat{u},0)$ such that $\hat{h}_{\hat{u}}(\hat{u},0)$
is small but not zero in absolute value, the inequality will hold.
QED Claim.

To be definite we choose the plus signs for the square roots.

Theorem 2, Section 2 (\ref{Theorem-2.-(}) gives a sufficient condition
for the embedding of a metric $\omega=Edu^{2}+Gdv^{2}$ where
\[
E(u,v)=\frac{1}{2}(a'(x+y))^{-2}.
\]
\[
G(u,v)-1=\frac{1}{2}(b'(x-y))^{-2}.
\]

In Theorem 2, $a(x)$ and $b(x)$ are arbitrary $C^{1}$ initial data
for two PDEs of first order. Therefore we may choose the initial data
$a(x),b(x)$ to be the solution functions of the initial value problem
for the ODE system. Therefore, by Theorem 2, $\omega$ thus defined
is isometrically embeddable. It remains to show that $\omega$ is
the given metric, i.e., we 

Claim: $E=R$ and $G=S.$ That is, 
\begin{equation}
\frac{1}{2}(a'(x+y))^{-2}=u_{\hat{v}}^{-2}\frac{\hat{G(}f(u,v),g(u,v))}{\hat{G(}f(u,v),g(u,v))+1}.\label{eq:E=00003DR}
\end{equation}
\begin{equation}
\frac{1}{2}(b'(x-y))^{-2}=\hat{h}_{\hat{v}}^{-2}(f(u,v),g(u,v))\frac{[\hat{G(}f(u,v),g(u,v))]^{2}}{\hat{G(}f(u,v),g(u,v))+1}-1.\label{G=00003DS}
\end{equation}
 Recall that functions $f$ and $g$ were defined in Section 3 depending
on solutions of initial value problems for PDEs (\ref{eq:linear PDE for u.},
\ref{linear PDE for v.}).

Now apply the parametric forms above of the solutions in Lemma 1,
Section 2. From the solutions to the PDEs of Lemma 1 in parametric
form we obtain for the left side of (\ref{eq:E=00003DR}) 
\[
\frac{1}{2}(a'(s))^{-2}.
\]
 However, by the solution of the ODE system above using our choice
of initial condition for the first PDE (\ref{eq:linear PDE for u.})$\,$$(i.e.,u=\hat{u}+\hat{v},$
cf. Remark Section 3) the left side is equal to

\[
\frac{\hat{G(}f(a(s),b(s)),g(a(s),b(s)))}{\hat{G(}f(a(s),b(s)),g(a(s),b(s))+1}.
\]
 Again using the parametric solutions of the PDEs of Lemma 1 ,i.e.,
$u=a(s),$ $v=b(s),$ we obtain for the left side of equation (\ref{eq:E=00003DR}):
\[
\frac{\hat{G(}f(u,v),g(u,v))}{\hat{G(}f(u,v),g(u,v)+1}.
\]

Which is the right side of (\ref{eq:E=00003DR}.) QED (\ref{eq:E=00003DR}).

Now apply the parametric forms of the PDEs in Lemma 1 to (\ref{G=00003DS}):
For the left side of (\ref{G=00003DS}) we obtain
\[
\frac{1}{2}(b'(\sigma)^{-2}.
\]
 In the ODE system the parameter $s$ may be replaced by $\mbox{}$$\sigma$.
The left side of (\ref{G=00003DS}) yields
\[
\hat{h}_{\hat{v}}^{-2}(f(a(\sigma),b(\sigma)),g(a(\sigma),b(\sigma)))\frac{[\hat{G(}f(a(\sigma),b(\sigma)),g(a(\sigma),b(\sigma)))]^{2}}{\hat{G(}f(a(\sigma),b(\sigma)),g(a(\sigma),b(\sigma)))+1}-1.
\]
 Since the equations for the solutions of the PDEs are the same as
the equations for the initial conditions we may use the conditions
given by $y=0.$ This yields 
\[
x=s=\sigma.
\]
 That is, on the initial plane $x$ is the parameter. Evaluate the
preceding expression on $y=0:$
\[
\hat{h}_{\hat{v}}^{-2}(f(a(s),b(s)),g(a(s),b(s)))\frac{[\hat{G(}f(a(s),b(s)),g(a(s),b(s)))]^{2}}{\hat{G(}f(a(s),b(s)),g(a(s),b(s)))+1}-1.
\]
 Again using the parametric solutions of the PDEs we obtain for the
left side
\[
\hat{h}_{\hat{v}}^{-2}(f(u,v),g(u,v))\frac{[\hat{G(}f(u,v),g(u,v))]^{2}}{\hat{G(}f(u,v),g(u,v))+1}-1.
\]
 which is the right side of (\ref{G=00003DS}). QED claim.

This concludes the proof of the isometric embedding theorem.

\section{An example to illustrate the method.}

Given the metric 
\[
\omega=d\hat{u}^{2}+d\hat{v}^{2},
\]
 i.e., the Euclidean/Riemannian metric, with $\hat{G}=1$ (cf. Section
3) . The first order PDEs become
\begin{equation}
u_{\hat{u}}-u_{\hat{v}}=0,\label{eq:sec3 pde for u}
\end{equation}
\begin{equation}
v_{\hat{u}}+v_{\hat{v}}=0.\label{eq:sec3pde for v}
\end{equation}
 A solution of (\ref{eq:sec3 pde for u}) is
\[
u=\hat{u}+\hat{v}
\]
 with initial value $u=\hat{u}$ at $\hat{v}=0.$ 

A solution of (\ref{eq:sec3pde for v}) is 
\[
v=\epsilon(\hat{u}-\hat{v})
\]
 with initial value $v=\epsilon\hat{u}$ at $\hat{v}=0.$

These two solutions give a transformation of parameters from the given
parameter system $\hat{u},\hat{v}$ to a new parameter system $u,v$
for $\omega.$ The Jacobian of the transformation has the value $J$=-2$\epsilon.$
We choose $\epsilon\neq0.$ Then the inverse parameter transformation
exists locally, is one-to-one and is in differentiability class $C^{1}.$
The inverse transformation can be found explicitly; for most examples
the explicit inverse cannot be found, we only know that it exists.
The inverse is
\[
\hat{u}=\frac{u+v}{2\epsilon},\,\hat{v}=\frac{u-v}{2\epsilon}.
\]
By definitions (\ref{eq:Def of f},\ref{eq:Def of g}) 
\[
f(u,v)=\frac{u+v}{2\epsilon},\, g(u,v)=\frac{u-v}{2\epsilon}.
\]
 Now apply this to the ODE system using $u=a(s),\, v=b(s),\,\hat{G}=1:$
\[
\frac{1}{2}(a'(s))^{-2}=\frac{1}{2}.
\]
\[
\frac{1}{2}(b'(s))^{-2}+1=[\hat{h}_{\hat{v}}(f(a(s),b(s)),g(a(s),b(s)))]^{-2}\frac{1}{2}.
\]
 Or
\[
[a'(s)]^{2}=1.
\]
\[
[b'(s)]^{-2}+2=[\hat{h}_{\hat{v}}(\frac{a(s)+b(s)}{2\epsilon},\frac{a(s)-b(s)}{2\epsilon})]^{-2}.
\]
 Put the system in the standard form $a'=...,\, b'=...$ choosing
the plus sign for the square roots for uniformity. Then
\[
a'(s)=1.
\]
\[
b'(s)=\left([\hat{h}_{\hat{v}}(\frac{a(s)+b(s)}{2\epsilon},\frac{a(s)-b(s)}{2\epsilon})]^{-2}-2\right)^{-\frac{1}{2}}.
\]
 A solution of the first equation is $a(s)=s.$ Substitute into the
second equation:
\[
b'(s)=\left([\hat{h}_{\hat{v}}(\frac{s+b(s)}{2\epsilon},\frac{s-b(s)}{2\epsilon})]^{-2}-2\right)^{-\frac{1}{2}}.
\]
\[
\]
 Solve the initial value problem for the second equation with initial
point $s=0,$ initial value $b(0)=0.$ 

The function $\hat{h}_{\hat{v}}$ depends on the choice of the initial
condition $\hat{h}(\hat{u},0)$ for the PDE $v_{\hat{u}}+v_{\hat{v}}=0$
for which a solution is 
\[
\hat{h}(\hat{u},\hat{v})=\epsilon(\hat{u}-\hat{v}).
\]
Therefore, for all $\hat{u},\hat{v}$, 
\[
\hat{h}_{\hat{v}}(\hat{u},\hat{v})=-\epsilon.
\]
 Thus the second equation becomes
\[
b'(s)=([(-\epsilon)]^{-2}-2)^{-\frac{1}{2}}=(\epsilon^{-2}-2)^{-\frac{1}{2}}
\]
which is real if $\left|\epsilon\right|$ is small. Thus a solution
of the ODE system is
\[
a(s)=s,
\]
\[
b(s)=(\epsilon^{-2}-2)^{-\frac{1}{2}}s.
\]
By the parametric solutions of the PDEs at the beginning of the preceding
section 4, 
\[
x+y=s,
\]
\[
x-y=\sigma,
\]
\[
u=s,
\]
\[
v=(\epsilon^{-2}-2)^{-\frac{1}{2}}\sigma.
\]
 Eliminating $s$ and $\sigma$ we find a relation between the $x,y$
and $u,v$ parameters: 
\[
u=x+y,
\]
\[
v=(\epsilon^{-2}-2)^{-\frac{1}{2}}(x-y).
\]
 The inverse tansformation is 
\[
x=H(u,v)=\frac{u((\epsilon^{-2}-2)^{-\frac{1}{2}}+v}{2(\epsilon^{-2}-2)^{-\frac{1}{2}}},
\]
\[
y=H^{*}(u,v)=\frac{v-u(\epsilon^{-2}-2)^{-\frac{1}{2}}}{-2(\epsilon^{-2}-2)^{-\frac{1}{2}}}.
\]
 By Theorem 2, the isometric embedding of $\omega$ in level parameters
$u,v$ is
\[
X(u,v)=\frac{u((\epsilon^{-2}-2)^{-\frac{1}{2}}+v}{2(\epsilon^{-2}-2)^{-\frac{1}{2}}}\mathbf{\mathrm{\mathbf{i\mathrm{}}+\mathbf{\mathrm{\mbox{\ensuremath{\frac{\mathrm{\mbox{\mbox{v}}-\mbox{u}(\epsilon^{-2}-2)^{-\frac{1}{2}}}}{-2(\epsilon^{-2}-2)^{-\frac{1}{2}}}}}\mathbf{j}+\mbox{\ensuremath{v}}\mathbf{k.}}}}}
\]

Referring now to section 2 we find 
\[
E(u,v)=\frac{1}{2}(a'(s))^{-2}=\frac{1}{2}.
\]
\[
G(u,v)=1+\frac{1}{2}(b'(\sigma))^{-2}=1+\frac{1}{2}((\epsilon^{-2}-2)^{-\frac{1}{2}})^{-2}=1+\frac{1}{2}(\epsilon^{-2}-2).
\]
 This may be confirmed by direct calculation from our formula for
$X(u,v).$ The metric in level parameters is
\[
\omega=\frac{1}{2}du^{2}+\frac{1}{2}[\epsilon^{-2}]dv^{2}.
\]
 Acknowledgements.

The author thanks the following for their advice and criticism:

Qing Han, Robert Bryant, Marcus Khuri, with special gratitude to Louis
Nirenberg and Jozef Dodziuk.

\[
\]

\[
\]
 
\[
\]
 
\[
\]
\[
\]
\[
\]

\[
\]
\[
\]

\title{
\[
\]
\[
\]
}

\[
\]

\[
\]

\end{document}